\documentclass[12pt]{article}
\usepackage{amsmath,graphicx,amssymb,amsthm}
\allowdisplaybreaks
\def\BB{\mathfrak B}

\def\NN{\mathfrak N}
\def\GN{\mathfrak{GN}}
\def\A{\mathcal A}
\def\B{\mathcal B}

\def\K{\mathcal K}

\def\B{\mathcal B}

\def\H{\mathcal H}
\def\K{\mathcal K}
\def\L{\mathcal L}
\def\M{\mathcal M}
\def\N{\mathcal N}

\def\S{\mathcal S}
\def\T{\mathcal T}

\DeclareMathOperator{\Aut}{Aut}
\def\amslatex{$\mathcal{A}\kern-.1667em\lower.5ex\hbox{$\mathcal{M}$}\kern-.125em\mathcal{S}$-\LaTeX}

\def\tensor{\mathop{\bar\otimes}}
\newtheorem{set}{set}[section]
\newtheorem{Corollary}[set]{Corollary}

\newtheorem{Lemma}[set]{Lemma}

\newtheorem{Proposition}[set]{Proposition}
\newtheorem{Remark}[set]{Remark}
\newtheorem{Theorem}[set]{Theorem}
\newcommand{\define}{\mathrel{\hbox{$\equiv$\hskip -.90em \lower .47ex \hbox{$\leftharpoondown$}}}}
\newcommand{\enifed}{\mathrel{\hbox{$\equiv$\hskip -.90em \lower .47ex \hbox{$\rightharpoondown$}}}}

\date{}
\begin{document}
\title {On Completely Singular von Neumann Subalgebras}

%\vskip1.3cm

\author{Junsheng Fang\sc \footnote{jfang@cisunix.unh.edu} }

%\author{\sc \footnote{Corresponding Author} \and \sc \and \sc }

\maketitle

\pagestyle{myheadings}

\begin{abstract}Let $\M$ be a
von Neumann algebra acting on a Hilbert space $\H$, and $\N$ be a
singular von Neumann subalgebra of $\M.$ If $\N\tensor\B(\K)$ is
singular in $\M\tensor\B(\K)$ for every Hilbert space $\K$, $\N$ is
said to be \emph{completely singular} in $\M$. We prove that if $\N$
is a singular abelian von Neumann subalgebra or if $\N$ is a
singular subfactor of a type ${\rm II}_1$ factor $\M$, then $\N$ is
completely singular in $\M$.  If $\H$ is separable, we prove that
$\N$ is completely singular in $\M$ if and only if for every
$\theta\in \Aut(\N')$ such that $\theta(X)=X$ for all $X\in\M'$,
then $\theta(Y)=Y$ for all $Y\in\N'$. As the first application, we
prove that if $\M$ is separable (with separable predual) and $\N$ is
completely singular in $\M$, then $\N\tensor\L$ is completely
singular in $\M\tensor \L$ for every separable von Neumann algebra
$\L$. As the second application, we prove that if $\N_1$ is a
singular subfactor of a type ${\rm II}_1$ factor  $\M_1$ and $\N_2$
is a completely singular von Neumann subalgebra of $\M_2$, then
$\N_1\tensor \N_2$ is completely singular in $\M_1\tensor \M_2$.
\end{abstract}

\vskip1.0cm

\section{Introduction}
 Let $\M$ be a von Neumann algebra acting on a Hilbert space $\H$.
A von Neumann subalgebra $\N$ of $\M$ is \emph{singular} if the only
unitary operators in $\M$ satisfying the condition $U\N U^*=\N$ are
those in $\N$. The study of singular von Neumann subalgebras has a
long and rich history (see for instance \cite{Di,J-P,Po1,Po2,Pu}).
Recently, there is a remarkable progress on  singular MASAs (maximal
abelian von Neumann subalgebras) in type ${\rm II}_1$ factors (see
\cite{S-S,RSS,SW}). In~\cite{S-S}, Allan Sinclair and Roger Smith
introduced a concept of asymptotic homomorphism property.
In~\cite{RSS}, a concept of weak asymptotic homomorphism property is
introduced. Let $\M$ be a type ${\rm II}_1$ factor and $\N$ be a von
Neumann subalgebra of $\M$. Then $\N\subseteq \M$ has the weak
asymptotic homomorphism property if for all $X_1,\cdots,X_n\in\M$
and $\epsilon>0$, there exists a unitary operator $U\in \N$ such
that
\[\|E_{\N}(X_iUX_j^*)-E_{\N}(E_\N(X_i)UE_\N(X_j)^*)\|_2<\epsilon.
\] Remarkably, in~\cite{SW}, it was shown that every singular MASA
 in a type ${\rm II}_1$ factor satisfies the weak asymptotic
homomorphism property. As a corollary, the tensor product of
singular MASAs in type ${\rm II}_1$ factors is proved to be a
singular MASA in the tensor product of type ${\rm II}_1$ factors
(see~\cite{SW}), which is a
well-known  hard question for long time. \\

It is very natural to ask the following question:
 if $\N_1$ and $\N_2$ are singular von
Neumann subalgebras of $\M_1$ and $\M_2$, respectively,  is
$\N_1\tensor \N_2$ singular in $\M_1\tensor\M_2$?  It turns out this
is not always true.  Let $\M_1=M_3(\mathbb{C})$ and
$\N_1=M_2(\mathbb{C})\oplus \mathbb{C}$.  Then
$P=\left(\begin{array}{ccc} 1&0&0\\
0&1&0\\
0&0&0\end{array}\right), Q=\left(\begin{array}{ccc} 0&0&0\\
0&0&0\\
0&0&1\end{array}\right)$ are central projections in $\N_1$ and
$\N_1=\{P,Q\}'$. Suppose $U\in\M_1$ is a unitary matrix such that
$U\N_1 U^*=\N_1$. Then $UPU^*=P$ and $UQU^*=Q$ (because $ad U$
preserves the center of $\N_1$ and $\tau(P)=\frac{2}{3}$,
$\tau(Q)=\frac{1}{3}$, where $\tau$ is the normalized trace on
$M_3(\mathbb{C})$). So $U\in \{P,Q\}'=\N_1$. This implies that
$\N_1$ is singular in $\M_1$. Let $\M_2=\B(l^2(\mathbb{N}))$ and
$\N_2=\M_2$. Then $\N_1\tensor
\B(l^2(\mathbb{N}))=M_2(\mathbb{C})\tensor\B(l^2(\mathbb{N}))\oplus
\mathbb{C}\tensor \B(l^2(\mathbb{N}))$ is not singular in
$\M_1\tensor\B(l^2(\mathbb{N}))=M_3(\mathbb{C})\tensor
\B(l^2(\mathbb{N}))$. Indeed, let $V$ be an isometry from
$l^2(\mathbb{N})$ onto $\mathbb{C}^2\otimes l^2(\mathbb{N})$, then
$U=\left(\begin{array}{cc} 0&V\\
V^*&0
\end{array}\right)$ is a unitary operator in $M_3(\mathbb{C})\tensor
 \B(l^2(\mathbb{N}))$ such that
$U(\N_1\tensor\B(l^2(\mathbb{N})))U^*=\N_1\tensor\B(l^2(\mathbb{N}))$.
Since $U$ is not in $\N_1\tensor\B(l^2(\mathbb{N}))$,
$\N_1\tensor\B(l^2(\mathbb{N}))$ is not singular in $\M_1\tensor
\B(l^2(\mathbb{N}))$. Indeed, $\N_1\tensor\B(l^2(\mathbb{N}))$ is
regular in $\M_1\tensor \B(l^2(\mathbb{N}))$ (see Remark 2.14). \\

Let $\M$ be a von Neumann algebra and $\N$ be a von Neumann
subalgebra of $\M.$ If $\N\tensor\B(\K)$ is singular in
$\M\tensor\B(\K)$ for every Hilbert space $\K$, then  $\N$ is said
to be \emph{completely singular} in $\M$. In section 2, we prove
that if $\N$ is a singular MASA  or if $\N$ a singular subfactor of
a type ${\rm II}_1$ factor $\M$, then $\N$ is completely singular in
$\M$. For every type ${\rm II}_1$ factor $\M$, we construct a
singular von Neumann subalgebra $\N$ of $\M$ ($\N\neq \M$) such that
$\N\tensor\B(l^2(\mathbb{N}))$ is regular  in $\M\tensor
\B(l^2(\mathbb{N}))$. Motivated by Lemma 1.2 of~\cite{Fa}, we obtain
a nice characterization of complete singularity in section 3. As the
first application, in section 4.1, we prove that  if $\M$ is
separable and $\N$ is completely singular in $\M$, then
$\N\tensor\L$ is completely singular in $\M\tensor \L$ for every
separable von Neumann algebra $\L$. As the second application, we
prove that if $\N_1$ is a singular subfactor of a type ${\rm II}_1$
factor
 $\M_1$ and $\N_2$ is a completely singular von Neumann subalgebra
of $\M_2$, then $\N_1\tensor \N_2$ is singular in $\M_1\tensor
\M_2$. The following question seems to be interesting: if $\N_1$,
$\N_2$ are completely singular von Neumann subalgebras of $\M_1$ and
$\M_2$, is $\N_1\tensor \N_2$ completely singular in
$\M_1\tensor\M_2$?

\section{On singularity and complete singularity}
\subsection{Normalizer and Normalizing groupoid of $\N$ in $\M$} Let $\M$ be a von
Neumann algebra and $\N$ be a von Neumann subalgebra of $\M$.  Then
$\NN_\M(\N)$ denotes the normalizer of $\N$ in $\M$:
\[\NN_\M(\N)=\{U\in\M:\,\, U\,\,  \text{is a unitary operator}, U\N U^*=\N\},\]
and $\GN_\M(\N)$ denotes the normalizing groupoid of $\N$ in $\M$:
\[\GN_\M(\N)=\{V\in\M:\,\, V\, \text{is a partial isometry with initial space
}\, E \,\text{and final space}\, F\]\[  \text{ such that}\, E,
F\in\N\, \text{and}\, V\N_E V^*=\N_F\},\] where $\N_E=E\N E$ and
$\N_F= F\N F$. $\N$ is singular in $\M$ if and only if
$\NN_\M(\N)''$, the von Neumann algebra generated by $\NN_\M(\N)$,
is $\N$. Recall that $\N$ is \emph{regular} in $\M$ if
$\NN_\M(\N)''=\M$.
\\

 If $\M$ is a finite von Neumann algebra and
$\N$ is a maximal abelian von Neumann subalgebra of $\M$, then $V\in
\GN_\M(\N)$ if and only if there is a unitary operator $U\in
\NN_\M(\N)$ and a projection $E\in\N$ such that $V=UE$ (\cite{J-P},
Theorem 2.1). In other words: any partial isometry that normalizes
$\N$ extends to a unitary operator that normalizes $\N$. As a
corollary, we have $\GN_\M(\N)''=\NN_\M(\N)''$, i.e., the von
Neumann algebra generated by the normalizing groupoid of $\N$ in
$\M$ is the von Neumann algebra generated by the normalizer of $\N$
in $\M$. If $\M$ is an infinite factor, e.g., type ${\rm III}$, and
$\N=\M$, then there is an isometry in $\M$ which can not be extended
to a unitary operator in $\M$. The following example tells us that
even the weak form $\GN_\M(\N)''=\NN_\M(\N)''$ can fail. Let
$\M=M_3(\mathbb{C})$ and $\N=M_2(\mathbb{C})\oplus \mathbb{C}$. As
we point out in the introduction, $\N$ is singular in $\M$, i.e.,
$\NN_\M(\N)''=\N$. Simple computations show that
$V=\left(\begin{array}{ccc} 0&0&0\\
0&0&1\\
0&1&0
\end{array}\right)$ is in $\GN_\M(\N)$. Note that $V$ is not in
$\N$.\\

Let  $V_1, V_2\in\M$ be two partial isometries in $\GN_\M(\N)$ and
$E_i=V_i^*V_i\in\N$, $i=1,2$. We say $V_1\preceq V_2$ if $E_1\leq
E_2$ and $V_1=V_2E_1$. It is obvious that $\preceq$ is a partial
order on the set of partial isometries in $\GN_\M(\N)$. Let
$\{V_\alpha\}$ be a totally ordered subset of $\GN_\M(\N)$, then
$V=\lim_\alpha V_\alpha$ (in strongly operator topology) exists and
$V\in \GN_\M(\N)$.

\begin{Lemma} If $\M$ is a finite von Neumann algebra and $\N$ is a
subfactor of $\M$, then for every $V\in \GN_\M(\N)$, there is a
unitary operator $U\in \NN_\M(\N)$ such that $V\preceq U$. In
particular, $\GN_\M(\N)''=\NN_\M(\N)''$.
\end{Lemma}
\begin{proof} By Zorn's lemma, there is a maximal
element $W\in \GN_\M(\N)$ such that $V\preceq W$. Let $E=W^*W$ and
$F=WW^*$. Then $E, F\neq 0$ and $E, F\in\N$.  We need to prove
$E=I$. If $E\neq I$, then $F\neq I$ since $\M$ is finite. So $I-E,
I-F\in\N$ are not $0$. Since $\N$ is a factor, there is a partial
isometry $V_1\in \N$ with initial space $E_1$, a non-zero
subprojection of $I-E$, and final space $E_2$, a non-zero
subprojection of $E$. Let $F'$ be the range space of $WE_2$. Then
$F'=WE_2W^*\in\N$. Since $\N$ is a factor, there is a paritial
isometry $V_2\in \N$ with initial space $F_2$, a non-zero
subprojection of $F'$, and final space $F_1$, a non-zero
subprojection of $I-F$. Now $W'=V_2WV_1$ is a partial isometry with
initial space $E_1\leq I-E$ and final space $F_1\leq I-F$. Simple
computation shows that $W+W'\in \GN_\M(\N)$. Note that $V\preceq
W\preceq W+W'$ and $W\neq W+W'$. It contradicts to the maximality of
$W$.
\end{proof}

\begin{Lemma}Let $\M$ be a von Neumann algebra and $\N$ be an abelian
 von Neumann subalgebra of $\M$.
  Then $\GN_\M(\N)''=\NN_\M(\N)''$.
\end{Lemma}
\begin{proof}Let $\M_1=\NN_\M(\N)''$. We only need to prove that $
\GN_\M(\N)''\subseteq \M_1$. For $V\in\M$ a partial isometry, define
$\S(V)=\{W\in\M_1: W$ is a partial isometry and $W\preceq V\}$.
Suppose $V\notin\M_1$. By Zorn's lemma, we can choose a maximal
element $W\in \S(V)$ such that $V-W\neq 0$ and $\S(V-W)=\{0\}$.
Since $W\in\M_1$, $V\in\M_1$ if and only if $V-W\in\M_1$. Therefore,
we can assume that $V\neq 0$ and $\S(V)=\{0\}$. Let $E=V^*V$ and
$F=VV^*$. Then $E\neq 0$ and $F\neq
0$.\\

If $E=F$, let $U=V+(I-E)$. Then $U\in \NN_\M(\N)$ and $V=UE\in
\M_1$. It is a contradiction. If $E\neq F$, we can assume that
$E_1=E(I-F)\neq 0$ (otherwise consider $V^*$). Let $V_1=VE_1$ and
$F_1$ be the final space of $V_1$.  Then $V_1\in \GN_\M(\N)$ with
initial space $E_1\leq I-F$ and final space $F_1\leq F$ such that
$0\neq V_1\preceq V$. Let $U=V_1+V_1^*+(I-E_1-F_1)$. Then $U\in
\NN_\M(\N)$ and $V_1=UE_1\in \M_1$. Note that $V_1\neq 0$ and
$V_1\preceq V$. $\S(V)\neq \{0\}$. It is a contradiction.
\end{proof}

 If $\N$ is singular in $\M$ and $E\in\N$ is a projection,
  $\N_E$ ($=E\N E$) may be not singular in $\M_E$.
   For example, let $\M=M_3(\mathbb{C})$ and
$\N=M_2(\mathbb{C})\oplus \mathbb{C}$ and
\[E=\left(\begin{array}{ccc} 0&0&0\\
0&1&0\\
0&0&1
\end{array}\right)\in\N.\] Then $\N_E$ is not singular in $\M_E$.
But we have the following result.

\begin{Lemma} Let $\N$ be a singular von Neumann subalgebra of $\M$
and $E\in\N$ be a projection. If $\N$ is a countably decomposable,
properly infinite von Neumann algebra, then $\N_E$ is singular in
$\M_E$.
\end{Lemma}
\begin{proof}   Let $P$ be the central support of $E$ relative to
$\N$. Then there are central projections $P_1$, $P_2$ of $\N$ such
that $P_1+P_2=P$ and $P_1E$ is finite and $P_2E$ is properly
infinite. Let $E_1=P_1E$ and $E_2=P_2E$. Then the central supports
of $E_1$ and $E_2$ are $P_1$ and $P_2$, respectively. Since $P_1$ is
a properly infinite countably decomposable projection and $E_1$ is a
finite projection in $\N_{P_1}$ and the central support of $E_1$ is
$P_1$, $P_1$ is a countably infinite sum of projections $\{E_{1n}\}$
in $\N$, each $E_{1n}$ is equivalent to $E_1$ in $\N_{P_1}$(see for
instance, Corollary 6.3.12 of~\cite{K-R}, volume 2). For $n\in
\mathbb{N}$, let $W_{1n}$ be a partial isometry in $\N_{P_1}$ such
that $W_{1n}^*W_{1n}=E_{1n}$ and $W_{1n}W_{1n}^*=E_1$. Since $P_2$
and $E_2$ are properly infinite projections in $\N_{P_2}$ with same
central support $P_2$ and $\N_{P_2}$ is countably decomposable,
$P_2$ and $E_2$ are equivalent in $\N_{P_2}$ (see for instance,
Corollary 6.3.5 of~\cite{K-R}, volume 2). Since $P_2$ is properly
infinite in $\N$, it can be decomposed into a countably infinite sum
of projections $\{E_{2n}\}$, each $E_{2n}$ is equivalent to $P_2$
and hence to $E_2$. For $n\in \mathbb{N}$, let $W_{2n}$ be a partial
isometry in $\N_{P_2}$ such that $W_{2n}^*W_{2n}=E_{2n}$ and
$W_{2n}W_{2n}^*=E_2$. Let $W_n=W_{1n}+W_{2n}\in\N$. Then
$W_n^*W_n=E_{1n}+E_{2n}$ and $W_nW_n^*=E_1+E_2=E$. \\

Suppose $V$ is a unitary operator in $\M_E$ such that $V\N_E
V^*=\N_E$. Define $U=\sum_{n=1}^\infty W_n^*VW_n+(I-P_1-P_2)$. Then
$U$ is a unitary operator and $U^*=\sum_{n=1}^\infty
W_n^*V^*W_n+(I-P_1-P_2)$. For any $T\in \N$,
$UTU^*=\sum_{m,n=1}^\infty W_m^*VW_mTW_n^*V^*W_n +(I-P_1-P_2)T$.
Note that $W_mTW_n^*\in \N_E$, $VW_mTW_n^*V^*\in\N_E$. So
$UTU^*\in\N$. Similarly, $U^*TU\in\N$. Thus $U\in \NN_\M(\N)$. Since
$\N$ is singular in $\M$, $U\in \N$. Therefore,
$W_1^*VW_1=U(E_{1n}+E_{2n})\in \N$. So $V=EVE=W_1W_1^*VW_1W_1^*\in
\N_E$. This implies that $\N_E$ is singular in $\M_E$.
\end{proof}
\subsection{Singular MASA  and singular subfactor (of type
${\rm II}_1$ factor) are completely singular}
\begin{Theorem} Let $\N$ be a von Neumann subalgebra of $\M$ and
$\K$ be a Hilbert space. If $\GN_\M(\N)''=\NN_\M(\N)''$, then
$\NN_{\M\tensor\B(\K)}(\N\tensor\B(\K))''=\NN_\M(\N)''\tensor\B(\K)$.
\end{Theorem}
Combine Theorem 2.4, Lemma 2.1 and Lemma 2.2, we have the following
corollaries.
\begin{Corollary} If $\M$ is a type ${\rm II}_1$ factor and $\N$ is a
singular subfactor of $\M$, then $\N$ is completely singular in
$\M$.
\end{Corollary}

\begin{Corollary} If $\N$ is a singular MASA of a von Neumann
algebra $\M$, then $\N$ is completely singular in $\M$.
\end{Corollary}

To prove Theorem 2.4, we need the following lemmas. We consider
dim$\K=2$ first, which  motivates the general case.
\begin{Lemma}Let $U=\left(\begin{array}{cc}
A_1&A_2\\
A_3&A_4
\end{array}
\right)$ be a unitary operator in $\M\tensor M_2(\mathbb{C})$. Then
the following conditions are equivalent: \begin{enumerate}
\item
$U(\N\tensor M_2(\mathbb{C}))U^*=\N\tensor M_2(\mathbb{C})$;
\item $A_i X A_j^*\in \N$ and $A_i^* X A_j\in\N$ for all $X\in\N$,
$1\leq i, j\leq 4$.
\end{enumerate}
\end{Lemma}
\begin{proof} $U(\N\tensor M_2(\mathbb{C}))U^*=\N\tensor
M_2(\mathbb{C})$ if and only if $U(\N\tensor
M_2(\mathbb{C}))U^*\subseteq\N\tensor M_2(\mathbb{C})$ and
$U^*(\N\tensor M_2(\mathbb{C}))U\subseteq\N\tensor M_2(\mathbb{C})$.
$U(\N\tensor M_2(\mathbb{C}))U^*\subseteq\N\tensor M_2(\mathbb{C})$
if and only if $U\left(\begin{array}{cc}
X&0\\
0&0
\end{array}
\right)U^*, U\left(\begin{array}{cc}
0&X\\
0&0
\end{array}
\right)U^*, U\left(\begin{array}{cc}
0&0\\
X&0
\end{array}
\right)U^*, U\left(\begin{array}{cc}
0&0\\
0&X
\end{array}
\right)U^*\in \N$ for all $X\in\N$. Simple computations show that
$U(\N\tensor M_2(\mathbb{C}))U^*\subseteq\N\tensor M_2(\mathbb{C})$
if and only if $A_i X A_j^*\in \N$  for all $X\in\N$, $1\leq i,
j\leq 4$.
\end{proof}

Since the proof of the following lemma is similar to the proof of
Lemma 2.7, we omit the proof.

\begin{Lemma}Let $U=(A_{ij})$ be a unitary operator in $\M\tensor \B(\K)$. Then
the following conditions are equivalent: \begin{enumerate}
\item
$U(\N\tensor \B(\K))U^*=\N\tensor \B(\K)$;
\item $A_i X A_j^*\in \N$ and $A_i^* X A_j\in\N$ for all $X\in\N$,
$1\leq i, j\leq dim \K$.
\end{enumerate}
\end{Lemma}

Let $X=I$ and $i=j$ in 2 of Lemma 2.8. We have the following
corollary.

\begin{Corollary} Let  $U=(A_{ij})$ be
 a unitary operator in $\M\tensor \B(\K)$ such
that $U(\N\tensor \B(\K))U^*\\=\N\tensor \B(\K)$. If
$A_{ij}=V_{ij}H_{ij}$ is the polar decomposition of $A_{ij}$, then
$H_{ij}\in\N$, $1\leq i,j\leq dim \K$.
\end{Corollary}

\begin{Lemma} Let $\N$ be a von Neumann algebra and $H$ be a
positive operator in $\N$ and $E$ be the closure of the range space
of $H$. Then the strong-operator closure of $\T=\{HXH:\,\, X\in\N\}$
is $\N_E$($=E\N E$).
\end{Lemma}
\begin{proof} It is easy to see $\T\subseteq \N_E$.
Let $H=\int_{\mathbb{R}} \lambda d E(\lambda)$ and
$E_{n}=E([1/n,\infty))$.  Then $\lim_{n\rightarrow \infty}E_{n}=E$
in strong-operator topology. Set $H_{n}=E_{n}H+(I-E_{n})$. Then
$H_{n}$ is invertible in $\N$. For $X\in\N_E$, let
$X_n=H_{n}^{-1}(E_{n} X E_{n})H_{n}^{-1}\in \N$. Then
$HX_nH=HH_{n}^{-1}E_{n}XE_{n}H_{n}^{-1}H=E_{n}XE_{n}\rightarrow
EXE=X$ in strong-operator topology. Hence, the strong-operator
closure of $\T$ contains $\N_E$.
\end{proof}

\begin{Lemma} Suppose $\N$ is a von Neumann subalgebra of
$\M$ and $A\in \M$ satisfying $A\N A^*\subseteq \N$ and $A^*\N
A\subseteq \N$. Let $A=VH$ be the polar decomposition and $E=V^*V$,
$F=VV^*$. Then $H, E,F\in\N$ and $V\in \GN_\M(\N)$.
\end{Lemma}
\begin{proof} By the assumption, $A^* I A=H^2\in \N$. So $H\in
\N$ and $E=R(H)\in\N$, where $R(H)$ is the closure of range space of
$H$. By symmetry, $F\in \N$. Note that $AXA^*=VHXHV^*\subseteq F\N
F=\N_F$ for all $X\in\N$. By lemma 2.10, $V \N_E V^*\subseteq \N_F$.
By $A^*XA\subseteq \N$ for all $X\in \N$ and similar arguments, $V^*
\N_F V\subseteq \N_E$. So $ \N_F\subseteq V\N_E V^*$. Thus $V
\N_EV^*=\N_F$, i.e., $V\in \GN_\M(\N)$.
\end{proof}

\noindent\begin{proof}[The proof of Theorem 2.4] Let $U_1\in
\NN_\M(\N)$ and $V$ be a unitary operator in $\B(\K)$. Then
$U_1\otimes V\in \NN_{\M\tensor\B(\K)}(\N\tensor\B(\K))$. So
$\NN_{\M\tensor\B(\K)}(\N\tensor\B(\K))''\supseteq
\NN_\M(\N)''\tensor\B(\K)$.\\

Let $U=(A_{ij})$ be a unitary operator in $\M\tensor \B(\K)$ such
that $U (\N\tensor\B(\K))U^*=\N\tensor \B(\K)$. Let
$A_{ij}=V_{ij}H_{ij}$ be the polar decomposition of $A_{ij}$. By
Lemma 2.8, Corollary 2.9 and Lemma 2.11, $H_{ij}\in\N$ and
$V_{ij}\in \GN_\M(\N)$. By the assumption of Theorem 2.4, $V_{ij}\in
\NN_\M(\N)''$. So $U\in \NN_\M(\N)''\tensor\B(\K)$, i.e.,
$\NN_{\M\tensor\B(\K)}(\N\tensor\B(\K))''\subseteq
\NN_\M(\N)''\tensor\B(\K)$.
\end{proof}
\subsection{On singular but not completely singular von Neumann
subalgebras}
\begin{Proposition} If $\N$ is a singular but not a completely singular von Neumann subalgebra of
$\M$, then there is a von Neumann subalgebra $\M_1$ of $\M$ and a
Hilbert space $\K$ such that $\N \subsetneqq \M_1$, $\N$ is singular
in $\M_1$ and $\N\tensor\B(\K)$ is regular in $\M_1\tensor \B(\K)$.
\end{Proposition}
\begin{proof} Since $\N$ is not completely singular in $\M$, there
is a Hilbert space $\K$ such that
$\L=\NN_{\M\tensor\B(\K)}(\N\tensor \B(\K))''\supsetneqq \N\tensor
\B(\K)$. Since $\N\tensor \B(\K)\subseteq \L\subseteq \M\tensor
\B(\K)$,  $\L=\M_1\tensor \B(\K)$ for some von Neumann algebra
$\M_1$, $\N\subsetneqq \M_1\subseteq \M$. Since $\N$ is singular in
$\M$, $\N$ is singular in $\M_1$. Since $\M_1\tensor
\B(\K)=\NN_{\M\tensor\B(\K)}(\N\tensor \B(\K))$, $\N\tensor\B(\K)$
is regular in $\M_1\tensor \B(\K)$.
\end{proof}

\begin{Proposition} If $\M$ is a type ${\rm II}_1$ factor, then there is a
singular von Neumann subalgebra $\N$ of $\M$ such that $\N\neq \M$
and $\N\tensor\B(l^2(\mathbb{N}))$ is regular in $\M\tensor
\B(l^2(\mathbb{N}))$.
\end{Proposition}
\begin{proof}Let $\M_1$ be a type ${\rm I}_3$ subfactor of $\M$ and
$\M_2=\M_1'\cap\M$. Then $\M_2$ is a type ${\rm II}_1$ factor. We
can identify $\M$ with $M_3(\mathbb{C})\tensor\M_2$ and $\M_1$ with
$M_3(\mathbb{C})\tensor \mathbb{C} I$.  With this indentification,
let $\N=(M_2(\mathbb{C})\oplus \mathbb{C})\tensor\M_2$. Then
\[P=\left(\begin{array}{ccc} 1&0&0\\
0&1&0\\
0&0&0
\end{array}\right)\otimes I,\,\,\text{and}\,\, Q=\left(\begin{array}{ccc} 0&0&0\\
0&0&0\\
0&0&1
\end{array}\right)\otimes I\] are central projections in $\N$.
$\N=\{P,Q\}'\cap\M$ and $\{P,Q\}''$ is the center of $\N$.  Let
$U\in\M$ be a unitary operator such that $U\N U^*=\N$. Then
$U\{P,Q\}''U^*=\{P,Q\}''$. Let $\tau$ be the unique tracial state on
$\M$. Then $\tau(P)=\frac{2}{3}$ and $\tau(Q)=\frac{1}{3}$. So
$UPU^*=P$ and $UQU^*=Q$. This implies that $U\in \{P,Q\}'\cap\M=\N$
and $\N$ is singular in $\M$.\\

To see $\N\tensor \B(l^2(\mathbb{N}))$ is not singular in $\M\tensor
\B(l^2(\mathbb{N}))$, we identify $\M\tensor\B(l^2(\mathbb{N}))$
with $M_3(\mathbb{C})\tensor \B(l^2(\mathbb{N}))\tensor\M_2$ and
$\N\tensor \B(l^2(\mathbb{N}))$ with
$(M_2(\mathbb{C})\oplus\mathbb{C})\tensor
\B(l^2(\mathbb{N}))\tensor\M_2$. Let $V$ be an isometry from
$l^2(\mathbb{N})$ onto $\mathbb{C}^2\otimes l^2(\mathbb{N})$, then
$U=\left(\begin{array}{cc} 0&V\\
V^*&0
\end{array}\right)$ is a unitary operator in $M_3(\mathbb{C})\tensor
\B(l^2(\mathbb{N}))$ such that
$U((M_2(\mathbb{C})\oplus\mathbb{C})\tensor
\B(l^2(\mathbb{N})))U^*=(M_2(\mathbb{C})\oplus\mathbb{C})\tensor
\B(l^2(\mathbb{N}))$. So $U\otimes I$ is a unitary operator in the
normalizer of $\N\tensor\B(l^2(\mathbb{N}))$ but $U\otimes I\notin\N\tensor\B(l^2(\mathbb{N}))$.\\

By Proposotion 2.12, there is  a von Neumann subalgebra $\L$ of $\M$
such that $\N \subsetneqq \L$ and $\N\tensor\B(l^2(\mathbb{N}))$ is
regular in $\L\tensor \B(l^2(\mathbb{N}))$. Since
$(M_2(\mathbb{C})\oplus \mathbb{C})\tensor\M_2\subsetneqq
\L\subseteq M_3(\mathbb{C})\tensor\M_2$, by Ge-Kadison's splitting
theorem (see~\cite{G-K}), $\L=\L_1\tensor \M_2$ for some von Neumann
algebra $\L_1$ such that $M_2(\mathbb{C})\oplus
\mathbb{C}\subsetneqq \L_1\subseteq M_3(\mathbb{C})$. Since
$M_3(\mathbb{C})$ is the unique von Neumann subalgebra of
$M_3(\mathbb{C})$ satisfies the above condition,
$\L_1=M_3(\mathbb{C})$. This implies that $\N\tensor
\B(l^2(\mathbb{N}))$ is regular in $\M\tensor \B(l^2(\mathbb{N}))$.
\end{proof}

\begin{Remark}\emph{By the proof
of Proposition 2.13, $(M_2(\mathbb{C})\oplus \mathbb{C})\tensor
\B(l^2(\mathbb{N}))$ is regular in $M_3(\mathbb{C})\tensor
\B(l^2(\mathbb{N}))$.}
\end{Remark}
\begin{Remark}\emph{ Let $\M$ be a type ${\rm II}_1$ factor and $\N$ be the
singular von Neumann subalgebra constructed as in the proof of
Proposition 2.13. It is easy to see that $\N\tensor \N$ is not
singular in $\M\tensor\M$.}
\end{Remark}

\section{A characterization of complete singularity}

\begin{Theorem} Let $\M$ be a von Neumann algebra acting on a
separable Hilbert space $\H$ and $\N$ be a von Neumann subalgebra of
$\M$. Then the following conditions are equivalent.
\begin{enumerate}
\item $\N$ is completely singular in $\M$;
\item $\N\tensor\B(l^2(\mathbb{N}))$ is singular in $\M\tensor
\B(l^2(\mathbb{N}))$;
\item If $\theta\in \Aut(\N')$ and $\theta(X)=X$ for all $X\in \M'$,
then $\theta(Y)=Y$ for all $Y\in\N'$.
\end{enumerate}
\end{Theorem}
\begin{proof} $``3\Rightarrow 1"$. Let $\K$ be a Hilbert space and $U\in \M\tensor
\B(\K)$ be a unitary operator such that
$U(\N\tensor\B(\K))U^*=\N\tensor\B(\K)$. Note that
$(\M\tensor\B(\K))'=\M'\tensor \mathbb{C}I_{\K}$ and
$(\N\tensor\B(\K))'=\N'\tensor \mathbb{C}I_{\K}$.  $\theta=ad U\in
\Aut (\N'\tensor \mathbb{C}I_\K)$. Since $U\in \M\tensor\B(\K)$,
$\theta(X\otimes I_\K)=X\otimes I_\K$ for all $X\in\M'$. By the
assumption of 3, $Y\otimes I_\K=\theta(Y\otimes I_\K)=U(Y\otimes
I_\K)U^*$ for all $Y\in \N'\tensor \mathbb{C}I_\K$. This implies
that $U\in \N\tensor\B(\K)$. Therefore, $\N\tensor\B(\K)$ is
singular in $\M\tensor\B(\K)$.\\

$``1\Rightarrow 2"$ is trivial.\\

$``2\Rightarrow 3"$. By~\cite{Ha}, there is a separable Hilbert
space $\H_1$ and a faithful normal representation $\phi$ of $\N'$
such that $\phi(\N')$ acts on  $\H_1$ in standard form.  Let
$\theta_1=\phi\cdot\theta\cdot \phi^{-1}$. Then $\theta_1\in \Aut
\phi(\N')$ and $\theta_1(\phi(X))=\phi(X)$ for all $X\in \M'$.  Now
there is a unitary operator $U_1\in\B(\H_1)$ such that
$\theta_1(\phi(Y))=U_1\phi(Y)U_1^*$ for all $Y\in\N'$. Let $\N_1$
and $\M_1$ be the commutant of $\phi(\N')$ and $\phi(\M')$ relative
to $\H_1$, respectively. Then $\N_1$ is a von Neumann subalgebra of
$\M_1$. Since $\theta_1(\phi(X))=U_1\phi(X)U_1^*=\phi(X)$ for all
$X\in \M'$, $U_1\in \M_1$. Since $\theta=ad U_1\in \Aut \phi(\N')$,
$\theta=ad U_1\in \Aut \N_1$. Now we only need to prove that $\N_1$
is a singular von Neumann subalgebra of $\M_1$. Then $U_1\in\N_1$
and $\theta_1(\phi(Y))=U_1\phi(Y)U_1^*=\phi(Y)$ for all $Y\in \N'$.
This
implies that $\theta(Y)=Y$ for all $Y\in\N'$.\\

By~\cite{Di} (Theorem 3, page 61),
$\phi=\phi_3\cdot\phi_2\cdot\phi_1$, where $\phi_1(\N')=\N'\tensor
\mathbb{C}I_\K$, $\K=l^2(\mathbb{N})$,
$\phi_2(\N'\tensor\mathbb{C}I_\K)=(\N'\tensor\mathbb{C}I_\K)E$,
$E\in (\N'\tensor\mathbb{C}I_\K)'=\N\tensor \B(\K)$ and $\phi_3$ is
a spacial isomorphism. We may assume that $\phi_3=id$. Then
$\N_1=E(\N\tensor \B(\K))E$ and $\M_1=E(\M\tensor \B(\K))E$, where
$E\in \N\tensor \B(\K)$. By 2, $\N\tensor\B(\K)$ is singular in
$\M\tensor\B(\K)$.  Note that $\N\tensor\B(\K)$ is a countably
decomposable, properly infinte von Neumann algebra. By Lemma 2.3,
$\N_1$ is singular in $\M_1$.
\end{proof}

Note that in the proof of $``3\Rightarrow 1"$ of Theorem 3.1, we do
not need the assumption that $\H$ is a separable Hilbert space.\\

 Let $\M$ be a von Neumann algebra. A von Neumann subalgebra $\BB$
of $\M$ is called \emph{maximal injective} if it is injective and if
it is maximal with respect to inclusion in the set of all injective
von Neumann subalgebras of $\M$ (see~\cite{Po2}).

\begin{Proposition} If $\BB$ is a maximal injective von Neumann
subalgebra of $\M$, then $\BB$ is completely singular in $\M$.
\end{Proposition}
\begin{proof} We can assume that $\M$ acts on $\H$ in standard form.
  Then $\BB'$ is a minimal injective von Neumann algebra extension
of $\M'$ (see~\cite{Fa}, 1.3). Let $\theta\in \Aut(\BB')$ such that
$\theta(X)=X$ for all $X\in\M'$. Then $\theta(Y)=Y$ for all $Y\in
\BB'$ by Lemma 1.2 of~\cite{Fa}. By Theorem 3.1, $\BB$ is completely
singular in $\M$.
\end{proof}

\section{Completely singular von Neumann subalgebras in tensor
products of von Neumann algebras}

\subsection{}
 The proof of the following lemma is similar
to the proof of Lemma 6.6 of~\cite{S-Z}

\begin{Lemma} Let $\M$ be a separable von Neumann algebra and
$\N$ be a  singular von Neumann subalgebra of $\M$. If $\A$ is an
abelian von Neumann algebra, then $\N\tensor\A$ is a singular von
Neumann subalgebra of $\M\tensor \A$.
\end{Lemma}
\begin{proof} We can assume that $\M$ acts on a separable Hilbert
space $\H$ in standard form and $\A$ is countably decomposable. Then
there is a *-isomorphism  from $\A$ onto $L^{\infty}(\Omega,\mu)$
with $\mu$ a probability Radon measure on some compact space
$\Omega$. To the *-isomorphism $\A\rightarrow L^{\infty}(\Omega,
\mu)$ corresponds canonically a *-isomorphism $\Phi$ from
$\B(\H)\tensor \A$ onto $L^\infty (\Omega,\mu;\B(\H))$. Note that
$\Phi(\M\tensor\A)(\omega)=\M$ and $\Phi(\N\tensor\A)(\omega)=\N$
for almost all $\omega\in\Omega$. Let $U\in \M\tensor\A$ be a
unitary operator such that $U(\N\tensor\A)U^*= \N\tensor\A$. Then
$\Phi(U)=U(\omega)$ such that $U(\omega)$ is a unitary operator in
$\M$ for almost all $\omega\in\Omega$. By $U(\N\tensor\A)U^*=
\N\tensor\A$, we have $U(\omega)\N U(\omega)^*=\N$ for almost all
 $\omega\in\Omega$. Since $\N$ is singular in $\M$, $U(\omega)\in\N$
 for almost all $\omega\in\Omega$. Hence $U\in\N\tensor \A$.
\end{proof}
Since for every Hilbert space $\K$, $\M\tensor\A\tensor\B(\K)$ is
canonically isomorphic to $\M\tensor\B(\K)\tensor\A$. We have the
following corollary.
\begin{Corollary} Let $\M$ be a separable von Neumann algebra and
$\N$ be a completely singular von Neumann subalgebra. If $\A$ is an
abelian von Neumann algebra, then $\N\tensor\A$ is a completely
singular von Neumann subalgebra of $\M\tensor \A$.
\end{Corollary}

\begin{Theorem} Let $\M$ be a separable von Neumann algebra and
$\N$ be a completely singular von Neumann subalgebra.  Then
$\N\tensor\L$ is  completely singular in $\M\tensor \L$ for every
separable von Neumann algebra $\L$.
\end{Theorem}

\begin{proof} We can assume that $\M$ and $\L$ act on separable Hilbert spaces $\H$
and $\K$ in standard form, respectively. Let $\theta$ be in $\Aut
(\N'\tensor \L')$ such that $\theta(X\otimes Z)=X\otimes Z$ for all
$X\in\M'$ and $Z\in \L'$. Let $\A$ be the center of $\L'$. Then
$(\mathbb{C}I_\H\tensor \L')'\cap (\N'\tensor \L')=(\B(\H)\tensor
\L)\cap (\N'\tensor \L')=(\B(\H)\cap\N')\tensor
(\L\cap\L')=\N'\tensor \A$. So for $T\in \N'\tensor \A$ and $Z\in
\L'$,  $T (I_\H\otimes Z)=(I_\H\otimes Z) T$ and
$\theta(T)\theta(I_\H\otimes Z)=\theta(I_\H\otimes Z)\theta(T)$.
Since $\theta(I_\H\otimes Z)=I_\H\otimes Z$, $\theta(T)(I_\H\otimes
Z)=(I_\H\otimes Z)\theta(T)$. This implies that $\theta(T)\in
\N'\tensor \A$. So $\theta\in \Aut (\N'\tensor \A)$ when $\theta$ is
restricted on $\N'\tensor\A$ such that $\theta(X\otimes Z)=X\otimes
Z$ for all $X\in\M'$ and $Z\in\A$.\\

Consider the standard representation $\phi$ of $\A$ on a separable
Hilbert space $\K_1$. Then $\phi(\A)'=\phi(\A)$. By Corollary 4.2,
$\N\tensor \phi(\A)$ is completely singular in $\M\tensor \phi(\A)$.
On $\H\otimes \K_1$, $(\N\tensor \phi(\A))'=\N'\tensor \phi(\A)$ and
$(\M\tensor \phi(\A))'=\M'\tensor \phi(\A)$.
 Note that $\theta_1=(id\tensor\phi)\cdot\theta\cdot
 (id\tensor\phi^{-1})
 \in \Aut (\N'\tensor \phi(\A))$ and $\theta_1(X\otimes Z')=(id\tensor
 \phi)\cdot \theta(X\otimes \phi^{-1}(Z'))=(id\tensor \phi)(X\otimes
 \phi^{-1}(Z')=X\otimes Z'$ for all $X\in\M'$ and $Z'\in\phi(\A)$.
  By Theorem 3.1, $\theta_1(Y\otimes Z')=Y\otimes Z'$ for all
  $Y\in\M'$ and $Z'\in\phi(\A)$. This implies that
$\theta(Y\otimes \phi^{-1}(Z'))=Y\otimes \phi^{-1}(Z')$ for all
$Y\in\N'$ and $Z'\in\phi(\A)$. Let $Z'=I_{\K_1}$. Then
$\theta(Y\otimes I_{\K})=Y\otimes I_{\K}$ for all $Y\in \N'$. Hence
$\theta(Y\otimes Z)=Y\otimes Z$ for all $Y\in\N'$ and $Z\in \L'$. By
Theorem 3.1, $\N\tensor\L$ is completely singular in $\M\tensor \L$.
\end{proof}
\subsection{}
\begin{Theorem} Let $\M_i$ be a separable von Neumann algebra and
$\N_i$ be a completely singular von Neumann subalgbebra of $\M_i$,
$i=1,2$. If $\N_1$ is a factor, then $\N_1\tensor \N_2$ is
completely singular in $\M_1\tensor \M_2$.
\end{Theorem}
\begin{proof} We can assume that $\M_1$ and $\M_2$ act on separable
 Hilbert space $\H_1$ and $\H_2$ in standard form, respectively. Let $\theta$ be in $\Aut
(\N_1'\tensor \N_2')$ such that $\theta(X_1\otimes X_2)=X_1\otimes
X_2$ for all $X_1\in\M_1'$ and $X_2\in \M_2'$.\\

Since $\N_1$ is a singular subfactor in $\M_1$,
$\N_1'\cap\M_1=\N_1'\cap\N_1=\mathbb{C}I_{\H_1}$.  Note that
$(\M_1'\tensor \mathbb{C}I_{\H_2})'\cap
(\N_1'\tensor\N_2')=(\M_1\tensor\B(\H_2))\cap (\N_1'\tensor\N_2')
=(\M_1\cap \N_1')\tensor (\B(\H_2)\cap\N_2')=
\mathbb{C}I_{\H_1}\tensor \N_2'$. We have
$\theta(\mathbb{C}I_{\H_1}\tensor \N_2')= \theta((\M_1\cap
\N_1')\tensor (\B(\H_2)\cap\N_2'))=\theta((\M_1\tensor\B(\H_2))\cap
(\N_1'\tensor\N_2'))=\theta((\M_1'\tensor \mathbb{C}I_{\H_2})'\cap
(\N_1'\tensor\N_2'))= \theta(\M_1'\tensor \mathbb{C}I_{\H_2})'\cap
\theta(\N_1'\tensor\N_2')= (\M_1'\tensor \mathbb{C}I_{\H_2})'\cap
(\N_1'\tensor\N_2')=\mathbb{C}I_{\H_1}\tensor \N_2'$. Since $\N_2$
is completely singular in $\M_2$ and $\theta(I_{\H_1}\otimes
X_2)=I_{\H_1}\otimes X_2$ for all $X_2\in\M_2'$,
$\theta(I_{\H_1}\otimes Y_2)=I_{\H_1}\otimes Y_2$ for all
$Y_2\in\M_2'$ by Theorem 3.1. Therefore, $\theta(X_1\otimes
Y_2)=X_1\otimes Y_2$ for all $X_1\in\M_1'$ and $Y_2\in \N_2'$. By
Theorem 4.3, $\N_1\tensor \M_2$ is completely singular in
$\M_1\tensor \M_2$. Since $\theta(X_1\otimes Y_2)=X_1\otimes Y_2$
for all $X_1\in\M_1'$ and $Y_2\in \N_2'$, by Theorem 3.1,
$\theta(Y_1\otimes Y_2)=Y_1\otimes Y_2$ for all $Y_1\in\N_1'$ and
$Y_2\in \N_2'$. By Theorem 3.1 again, $\N_1\tensor \N_2$ is
completely singular in $\M_1\tensor \M_2$.
\end{proof}

Combining Theorem 4.4 and Corollary 2.5, we obtain the following
corollary, which generalizes Corollary 4.4 of~\cite{SWW}.
\begin{Corollary}If
$\N_1$ is a singular subfactor of a type ${\rm II}_1$ factor  $\M_1$
and $\N_2$ is a completely singular von Neumann subalgebra of
$\M_2$, then $\N_1\tensor \N_2$ is completely singular in
$\M_1\tensor \M_2$.
\end{Corollary}

\end{document}